\documentclass[11pt]{amsart}

\usepackage{tabularx}
\usepackage{amsmath}
\usepackage{amssymb}
\usepackage{amsfonts}
\usepackage{boxedminipage}
\usepackage{appendix}
\usepackage{algorithm2e}
\RestyleAlgo{boxruled}
\usepackage{mathtools}
\usepackage{physics}
\usepackage{subfigure}

 \usepackage[foot]{amsaddr}

\usepackage{cite}

\usepackage{graphicx}
\usepackage[export]{adjustbox}

\usepackage{flexisym}	 
\usepackage{breqn}       

\usepackage{amsmath,amssymb,amsthm}
\usepackage{mathtools}   % Some improvements to AMS-LaTeX
\usepackage{physics}

\usepackage[utf8]{inputenc} 

\usepackage{csquotes}     
\MakeOuterQuote{"}       

\usepackage{fullpage}

\usepackage[colorinlistoftodos]{todonotes}   

\usepackage{hyperref} % adds hyper links inside the generated pdf file
 \hypersetup{%
 	colorlinks=true,       % false: boxed links; true: colored links
 	linkcolor=blue,          % color of internal links
 	citecolor=blue,        % color of links to bibliography
 	filecolor=magenta,      % color of file links
 	urlcolor=blue         
 }
\usepackage{listings}
\usepackage{xcolor}

\definecolor{codegreen}{rgb}{0,0.6,0}
\definecolor{codegray}{rgb}{0.5,0.5,0.5}
\definecolor{codepurple}{rgb}{0.58,0,0.82}
\definecolor{backcolour}{rgb}{0.95,0.95,0.92}

\lstdefinestyle{mystyle}{
    backgroundcolor=\color{backcolour},   
    commentstyle=\color{codegreen},
    keywordstyle=\color{magenta},
    numberstyle=\tiny\color{codegray},
    stringstyle=\color{codepurple},
    basicstyle=\ttfamily\footnotesize,
    breakatwhitespace=false,         
    breaklines=true,                 
    captionpos=b,                    
    keepspaces=true,                 
    numbers=left,                    
    numbersep=5pt,                  
    showspaces=false,                
    showstringspaces=false,
    showtabs=false,                  
    tabsize=2
}

\usepackage{graphicx} % Added by Shoubhik Mondal, 2014.
\makeatletter
\renewcommand{\email}[2][]{%
  \ifx\emails\@empty\relax\else{\g@addto@macro\emails{,\space}}\fi%
  \@ifnotempty{#1}{\g@addto@macro\emails{\textrm{(#1)}\space}}%
  \g@addto@macro\emails{#2}%
}
\makeatother

\title{On the Role of Tikhonov Regularizations in Standard Optimization Problems}
\date{\today}
\author{J. Adriazola}
\address{Department of Mathematics, University of California, Santa Barbara, CA, USA, 93106}
\email{jadriazola@ucsb.edu}
\begin{document}
\maketitle
\begin{abstract}
Tikhonov regularization is a common technique used when solving poorly behaved optimization problems. Often, and with good reason, this technique is applied by practitioners in an ad hoc fashion. In this note, we systematically illustrate the role of Tikhonov regularizations in two simple, yet instructive examples. In one example, we use regular perturbation theory to predict the impact Tikhonov regularizations have on condition numbers of symmetric, positive semi-definite matrices. We then use a numerical example to confirm our result. In another example, we construct an exactly solvable optimal control problem that exhibits a boundary layer phenomenon. Since optimal control problems are rarely exactly solvable, this brings clarity to how vital Tikhonov regularizations are for the class of problems this example represents. We solve the problem numerically using MATLAB's built-in optimization software and compare results with a Galerkin-type/genetic algorithm. We find the second method generally outperforms MATLAB's optimization routines and we provide the MATLAB code used to generate the numerics. Finally, to demonstrate the utility of Tikhonov regularization from the second example in a more technologically realistic application, we show that regularization guides the design of gradient-index optical devices toward those which are significantly less expensive to fabricate.

\end{abstract}

\section{Introduction}\label{section:Intro}
Ill-posed and ill-conditioned problems frequently arise in applied and computational mathematics. To overcome this common challenge, one typically introduces a type of regularization that either converts the problem to a well-posed one or better conditions the resulting computational problem. The Tikhonov regularization technique is often chosen as the first strategy toward this goal and  has been successfully applied in several contexts ranging from ill-posed optimal control problems~\cite{Hintermuller,Hohenester,Mennemann,me,me2} to the numerical solution of poorly-conditioned integral equations~\cite{Tikhonov1995,Tik1,Tik2,Tik4,Tik5}. In statistics, the method is better known as ridge regression~\cite{Tik3}.

To illustrate how Tikhonov regularizations apply in a general optimization setting, consider the problem of finding $x_*$ such that
\begin{equation}
    J(x_*)=\min_{x\in X}J(x),
\end{equation}
subject to
\begin{equation}
    f(x)=0,
\end{equation}
where $J:X\to\mathbb{R}$ is called an objective functional, $f:X\to Y$ is a constraint function, and $X,\ Y$ are general Banach spaces. Although the problem may be solvable, it may either be too computationally intensive to solve or fail to admit solutions which remain feasible to implement in a scientific or technological application. Thus, one often seeks to define a modified problem
\begin{equation}\label{eq:regur}
    J_{\varepsilon}(x) = J(x) + \varepsilon K(x) 
\end{equation}
while still remaining faithful to the originally intended goal.

Of the many possible pitfalls that may be encountered when attempting to solve an ill-posed problem, one important example is where the optimal solution $x_*\in X$ lacks a desired degree of regularity or smoothness. A Tikhonov regularization enforces smoothness of solutions to this problem by modifying the objective to include a penalization. A typical choice is setting 
\begin{equation}\label{eq:H1}
    K(x)=\norm*{x}^2_{\dot{H}^1(\Omega)}
\end{equation} in Equation~\eqref{eq:regur}, where $\Omega\subset\mathbb{C}^n,\ $ $n\in\mathbb{Z}^+,$ and $\varepsilon\in\mathbb{R}^+.$ The notation $\dot{H}^1(\Omega)$ indicates a homogeneous Sobolev space with norm given by $\left\|D x\right\|_{L^{2}(\Omega)},$ where $L^{2}(\Omega)$ is the standard Lebesgue (Hilbert) space~\cite{Lieb} and $D$ indicates differentiation with respect to elements from $\Omega.$ Simply put, the introduction of the regularization penalizes solutions that rapidly vary over $\Omega.$ 

The problem of finding an optimal value of $\varepsilon$ is known in statistics as the so-called bias-variance dilemma~\cite{Neal}. Although an optimal choice of the Tikhonov parameter $\varepsilon$ is known in certain circumstances~\cite{Wahba}, the parameter is typically chosen in an ad hoc manner. Through the simple and instructive examples studied here, we reveal suitable regularization parameters in simpler settings which may guide choices in more complex ones. Of course, choosing suitable regularization parameters in more complex, real-world applications remains to be done on a case-by-case basis.

This work is organized as follows. Section~\ref{section:LA} shows how the Tikhonov technique better conditions a symmetric, positive semi-definite matrix whose inverse is desired. There, we quantitatively show the type of impact a small Tikhonov regularization can have on an extremely ill-conditioned matrix. In Section~\ref{section:OC}, we show an example of an optimal control problem admitting a closed-form solution. Since the problem is exactly solvable, this provides a testbed for standard optimization methods of which we compare two; an interior-point method~\cite{BoydV} implemented by MATLAB's \texttt{fmincon} and a method, due to Calarco, et al.~\cite{Caneva,Doria}, called the Chopped Random Basis (CRAB) method which we also implement in MATLAB. In Appendix~\ref{section:appendix}, we provide MATLAB code reflecting our implementation of the methods used to solve the optimal control problem discussed in Section ~\ref{section:OC}.
In Section~\ref{section:Beam}, we formulate and solve an optimal control problem that models the design of gradient-index optical devices. The problem is formulated, in the Schr\"odinger optics regime, as finding the optimal coupler between two waveguides with different transverse profiles. The problem is solved numerically and the utility of Tikhonov regularization is demonstrated in terms of the manufacturability of the resulting waveguides.

\section{Regularizing an Ill-Conditioned Matrix}\label{section:LA}

The first example we consider is matrix inversion of an ill-conditioned matrix. Let $n\in\mathbb{Z}^+$ and consider the following optimization problem:

\begin{equation}\label{eq:simpletik}
    \min_{x\in\mathbb{R}^n}J=\min_{x\in\mathbb{R}^n}\left\{\frac{1}{2}x^{\mathsf{T}}Ax-x^{\mathsf{T}}b+\frac{\varepsilon}{2}x^{\mathsf{T}}D^{\mathsf{T}}Dx\right\}
\end{equation}
where $A$ is some real, symmetric, and positive semi-definite $n\times n$ matrix, $b\in\mathbb{R}^n,$ $D$ is some, as of yet, unspecified $n\times m$ matrix, $m\in\mathbb{Z}^+,$ and $\varepsilon$ is a small, real number. The restrictions placed on $A$ are done so its spectrum can be analyzed in detail. 

By taking the gradient of the objective $J$ in Equation~\eqref{eq:simpletik} with respect to the variable $x$, we see that the optimal vector $x^*$ satisfies
\begin{equation}\label{eq:regeq}
    \left(A+\varepsilon D^{\mathsf{T}}D\right)x^*=b.
\end{equation}
Note that in the case $\varepsilon=0$, the solution to Equation~\eqref{eq:regeq} is given by
$x^*=A^{-1}b.$ If $A$ is singular, then the task of finding $x^*$ is ill-posed, and if $A$ is nearly-singular, computing $x^*$ is challenging and often requires substantial computational resources for large $n$. 

It is well-known that the condition number of a symmetric, positive semi-definite matrix $A$, with condition number denoted by $c(A),$ is given by the ratio of its largest eigenvalue to its smallest eigenvalue~\cite{StrangLA}. For this reason, we aim to better understand the effect of the Tikhonov regularizer $\varepsilon D^{\mathsf{T}}D$ on the spectrum of $A$ for small enough $\varepsilon$ since this informs us of the condition number of the regularized matrix $A+\varepsilon D^{\mathsf{T}}D$. 

To this end, further assume the eigenvalues of $A$ are simple, and consider the asymptotic expansions
\begin{subequations}
\begin{align}
    \lambda(\varepsilon)\sim\sum_{j=0}^{\infty}\lambda_j\varepsilon^j\qquad\mathrm{as}\ \varepsilon\to0,\\
    \nu(\varepsilon)\sim\sum_{j=0}^{\infty}\nu_j\varepsilon^j\qquad\mathrm{as}\ \varepsilon\to0,
\end{align}
\end{subequations}
where $\nu_0$ and $\lambda_0$ are defined as the respective eigenvector/eigenvalue pair of the matrix $A$. To leading order, we have 
\begin{equation}\label{eq:ground}
    A\nu_0=\lambda_0\nu_0,
\end{equation}
which is automatically satisfied by the definition of $\nu_0$ and $\lambda_0$. 

To next order, i.e., $\mathcal{O}(\varepsilon),$ we have
\begin{equation}\label{eq:nextorderTik}
    (A-\lambda_0I)\nu_1=(\lambda_1-D^{\mathsf{T}}D)\nu_0,
\end{equation}
where $I$ is the $n\times n$ identity matrix. The Fredholm alternative~\cite{StrangLA} requires the right-hand side to be orthogonal to the left-eigenvector of $A$ so that Equation~\eqref{eq:nextorderTik} is solvable. Since $A$ is symmetric, the left eigenvector of $A$ is simply $\nu_0^{\mathsf{T}}.$ The solvability condition of Equation~\eqref{eq:nextorderTik} is then given by
\begin{equation}\label{eq:fredcond}
    \lambda_1=\frac{\nu_0^{\mathsf{T}}D^{\mathsf{T}}D\nu_0}{\nu_0^{\mathsf{T}}\nu_0}.
\end{equation}
Thus, the resulting condition number is
%  \begin{equation}
%  \begin{alignedat}{2}
%  c(A+\varepsilon D^{\mathsf{T}}D)=\frac{\lambda^{\mathrm{max}}}{\lambda^{\mathrm{min}}}
%      &\sim\frac{\lambda_0^{\mathrm{max}}+\varepsilon\lambda_1^{\mathrm{max}}}{\lambda_0^{\mathrm{min}}+\varepsilon\lambda_1^{\mathrm{min}}}+\mathcal{O}\left(\varepsilon^2\right)&&\qquad\mathrm{as}\ \varepsilon\to0\\
%      &\sim\frac{\lambda_0^{\mathrm{max}}+\varepsilon\lambda_1^{\mathrm{max}}}{\lambda_0^{\mathrm{min}}}\left(1-\varepsilon\frac{\lambda_1^{\mathrm{min}}}{\lambda_0^{\mathrm{min}}}\right)+\mathcal{O}\left(\varepsilon^2\right)&&\qquad\mathrm{as}\ \varepsilon\to0\\
%      &\sim c(A)+\varepsilon\frac{\lambda_1^{\mathrm{max}}-c(A)\lambda_1^{\mathrm{min}}}{\lambda_0^{\mathrm{min}}}+\mathcal{O}\left(\varepsilon^2\right)&&\qquad\mathrm{as}\ \varepsilon\to0
%  \end{alignedat}
%  \end{equation}
\begin{equation}\label{eq:TikEffect}
 c(A+\varepsilon D^{\mathsf{T}}D)=\frac{\lambda^{\mathrm{max}}}{\lambda^{\mathrm{min}}}
     \sim\frac{\lambda_0^{\mathrm{max}}+\varepsilon\lambda_1^{\mathrm{max}}}{\lambda_0^{\mathrm{min}}+\varepsilon\lambda_1^{\mathrm{min}}}+\mathcal{O}(\varepsilon^2),\qquad\mathrm{as}\ \varepsilon\to0,
 \end{equation}
where $\lambda_0^{\rm min/max}$ are, respectively, the smallest and largest eigenvalues satisfying Equation~\eqref{eq:ground} while $\lambda_1^{\rm min/max}$ are the smallest and largest numbers satisfying Equation~\eqref{eq:fredcond}. Presumably, $\lambda_0^{\mathrm{min}}$ is nearly zero, so long as $\lambda_1^{\mathrm{min}}$ is not, the regularization has a high likelihood of being effective. 

Indeed, the asymptotic formula~\eqref{eq:TikEffect} gives a criteria for finding effective $A-$specific Tikhonov matrices, i.e., the Tikhonov matrix $D$ solves the following max-min problem
\begin{equation}\label{eq:maxmin}
    \max_{D\in\mathcal{D}}\min_{\nu_0\in\mathcal{V}}\left\{\frac{\nu_0^{\mathsf{T}}D^{\mathsf{T}}D\nu_0}{\nu_0^{\mathsf{T}}\nu_0}\right\},
\end{equation}
where $\mathcal{D}$ is the space of all $n\times m$ matrices and $\mathcal{V}$ is the set of eigenvectors of the ill-conditioned matrix $A.$ We do not take up a study of Equation~\eqref{eq:maxmin} here, and instead leave further investigation for future study.

To better illustrate the effect a Tikhonov regularization has on a nearly singular operator, take as an example the symmetric, positive-definite matrix
\begin{equation}
A=\begin{pmatrix}
    1 & 1\\
    1 & 1+\mu
\end{pmatrix},\qquad \mu>0
\end{equation}
which has eigenvalues 
\begin{equation}
    \lambda_0^{\rm min}=\frac{1}{2}\left(\mu +2-\sqrt{\mu ^2+4}\right),\quad \lambda_0^{\rm max}=\frac{1}{2}\left(\mu +2+\sqrt{\mu ^2+4}\right).
\end{equation}
 and corresponding eigenvectors
 \begin{equation}
 \nu_0^{\rm min}=\begin{pmatrix}
           \frac{1}{2}\left(-\mu-\sqrt{\mu ^2+4}\right)\\
           1
\end{pmatrix},\quad  \nu_0^{\rm max}=\begin{pmatrix}
           \frac{1}{2}\left(-\mu+\sqrt{\mu ^2+4}\right)\\
           1
\end{pmatrix},
 \end{equation}
Clearly, $A$ is ill-conditioned since $c(A)\to\infty$ as $\mu\to 0$ from above. Indeed, for $\mu=10^{-6},$ the condition number $c(A)= 4\times10^6$ to seven digits of precision. However, the asymptotic result~\eqref{eq:TikEffect} predicts that simply choosing $D$ to be the forward difference operator
\begin{equation}
    D=\begin{pmatrix}
    1 & 0\\
    -1 & 1\\
    0 & -1
\end{pmatrix},
\end{equation}
a choice consistent with Equation~\eqref{eq:H1} in this finite-dimensional setting, and a small regularization parameter of $\varepsilon=0.01,$ the regularized condition number $c(A+\varepsilon D^{\mathsf{T}}D)=203$ to three digits of precision; a noteworthy improvement. In fact, the asymptotic result is well within the same order of magnitude of the true value of the condition number 
\begin{equation}
    c(A+\varepsilon D^{\mathsf{T}}D)=\frac{2 + 4\varepsilon + \mu + \sqrt{4 - 8\varepsilon + 4\varepsilon^2 + \mu^2}}{2 + 4\varepsilon + \mu - \sqrt{4 - 8\varepsilon + 4\varepsilon^2 + \mu^2}}
\end{equation}
which is 67 to two digits of precision. 

Rounding out the example, we choose an example vector $b=\begin{pmatrix}
    \frac{1}{2},\ \frac{1}{2}
\end{pmatrix}^{\mathsf{T}}$ to demonstrate that we can still solve the matrix inversion in a meaningful way. When $\varepsilon=0$, the exact solution is given by
$x_{\rm unreg}=\begin{pmatrix}
    \frac{1}{2},\ 
    0
\end{pmatrix}^{\mathsf{T}},
$
while, to three digits of precision and when $\varepsilon=0.01,$ the solution
$x_{\rm reg}=\begin{pmatrix}
    0.248,
    \
    0.248
    \end{pmatrix}^{\mathsf{T}}.
$
The regularized solution has a residual, i.e., $J(x_{\rm reg})$ when $\varepsilon=0,$ of $1.24\times10^{-5}.$ The vast improvement on the condition number of $A,$ while still remaining faithful to the unregularized numerical computation, clearly demonstrates that even a small Tikhonov regularization can greatly mitigate an extremely ill-conditioned numerical computation.

% Choosing the operator $D$ is problem-specific. In practice, using second order centered-difference matrices $D^{\mathsf{T}}D$ is often sufficient. This is the choice we made implicitly in the following example

\section{Regularizing an Exactly Solvable Optimal Control Problem}\label{section:OC}

The second example we consider is an optimal control problem representative of optimization problems that frequently arise in controlling dispersive waves, e.g.~\cite{me,me2,Hohenester}. This optimal control problem, by virtue of its exact solution, helps reveal some of the typical numerical difficulties exhibited by this class of problems. Consider the following Lagrange control problem
\begin{equation}\label{eq:testcontprob}
    \min_{u\in\mathcal{U}}J=\min_{u\in\mathcal{U}}\int_0^{\pi}\left( x+\frac{1}{2}u^2+\frac{\varepsilon}{2}\left(\frac{\mathrm{d}u}{\mathrm{d}t}\right)^2\right)\mathrm{d}t,
\end{equation}
where the admissible class $\mathcal{U}=\left\{u\in H^1([0,\pi]):u(0)=u(\pi)=0\right\},$ subject to the forced harmonic oscillator
\begin{subequations}\label{eq:states}
\begin{alignat}{2}
    \frac{\mathrm{d}x}{\mathrm{d}t}&=p,\qquad &x(0)=0,\\
    \frac{\mathrm{d}p}{\mathrm{d}t}&=u-x,\qquad &p(0)=0.
\end{alignat}
\end{subequations}
The prescription of boundary conditions in defining the class of admissible controls $\mathcal{U}$ will soon play a critical role when considering the unregularized case, i.e., $\varepsilon=0$. For ease of notation in what follows, we use overhead dots to denote differentiation with respect to time, e.g., $\dot{x}=\frac{\mathrm{d}x}{\mathrm{d}t}.$

We recast the problem into a form more amenable to standard tools from the calculus of variations. Using Lagrange multipliers, or costates in the language of optimal control theory and denoted here by $\lambda(t)$ and $\mu(t)$, the equivalent unconstrained saddle-point problem is given by
\begin{align}
\begin{split}
&\min_{u\in\mathcal{U}}\int_0^{\pi}\mathcal{L}\left(x,\dot{x},p,\dot{p},\lambda,\mu,u,\dot{u}\right)\mathrm{d}t\\
=&\min_{u\in\mathcal{U}}\int_0^{\pi}\left(x+\frac{1}{2}u^2+\frac{\varepsilon}{2}\dot{u}^2+\lambda(\dot{x}-p)+\mu(\dot{p}+x-u)\right)\mathrm{d}t.
\end{split}
\end{align}
The Euler-Lagrange equations (see, for example,~\cite{Gelfand,Witelski,Bryson}) are given by
\begin{subequations}
\begin{alignat}{2}
\label{eq:ch2cos1}
\frac{\delta\mathcal{L}}{\delta x}&=\frac{\partial\mathcal{L}}{\partial{x}}-\frac{d}{dt}\frac{\partial\mathcal{L}}{\partial{\dot{x}}}=1-\dot{\lambda}+\mu=0,\qquad &\frac{\partial\mathcal{L}}{\partial{\dot{x}}}\bigg|_{t=\pi}=\lambda(\pi)=0,\\
\label{eq:ch2cos2}
\frac{\delta\mathcal{L}}{\delta p}&=\frac{\partial\mathcal{L}}{\partial{p}}-\frac{d}{dt}\frac{\partial\mathcal{L}}{\partial{\dot{p}}}=\lambda-\dot{\mu}=0,\qquad &\frac{\partial\mathcal{L}}{\partial{\dot{p}}}\bigg|_{t=\pi}=\mu(\pi)=0,\\
\label{eq:ch2state1}
\frac{\delta\mathcal{L}}{\delta \lambda}&=\frac{\partial\mathcal{L}}{\partial{\lambda}}=\dot{x}-p=0,\qquad &x(0)=0,\\
\label{eq:ch2state2}
\frac{\delta\mathcal{L}}{\delta \mu}&=\frac{\partial\mathcal{L}}{\partial{\mu}}=\dot{p}+x-u=0, \qquad &p(0)=0,\\
\label{eq:ch2control}
\frac{\delta\mathcal{L}}{\delta u}&=\frac{\partial\mathcal{L}}{\partial{u}}-\frac{d}{dt}\frac{\partial\mathcal{L}}{\partial{\dot{u}}}=u-\mu-\varepsilon\ddot{u}=0,\qquad &u(0)=u(\pi)=0.
\end{alignat}
\end{subequations}
Note that the costates' terminal conditions, present in Equations~\eqref{eq:ch2cos1} and~\eqref{eq:ch2cos2}, are a result of an integration by parts again test functions $\varphi\in C^\infty([0,T]),$ which do not necessarily vanish at $t=\pi,$ in the derivation of the Euler-Lagrange equations.
% More detail about the analogous derivation of these Euler-Lagrange equations in the context of PDE constraints is provided later on in this dissertation; for example in Appendix~\ref{section:StatOpt}.

By direct integration, it's easy to see the following functions
\begin{subequations}
\begin{align}
\lambda_*&=-\sin(t),\\
\mu_*&=-\cos(t)-1
\end{align}
\end{subequations}
are the optimal costates. This implies, after using the method of variation of parameters, see e.g.~\cite{Bender}, to solve Equation~\eqref{eq:ch2control}, the optimal control is
\begin{equation}\label{eq:closedcontrol}
u_*=\frac{\text{csch}\left(\frac{\pi }{\sqrt{\varepsilon }}\right) \left((\varepsilon +2) \sinh \left(\frac{\pi -t}{\sqrt{\varepsilon }}\right)+\varepsilon  \sinh \left(\frac{t}{\sqrt{\varepsilon }}\right)\right)-\cos (t)-\varepsilon-1}{\varepsilon +1}.
\end{equation}
Additionally, we solve for the optimal state $x_*$ in a similar fashion which results in
\begin{multline}\label{eq:closedstate}
x_*=\frac{\sqrt{\varepsilon } \text{csch}\left(\frac{\pi }{\varepsilon }\right)}{(\varepsilon +1)^2}\left(\varepsilon ^{3/2} \sinh \left(\frac{t}{\sqrt{\varepsilon }}\right)-\varepsilon  \sin (t)\right)-\frac{t \sin (t)}{2 (\varepsilon +1)}+\cos (t)-1\\
+\frac{\sqrt{\varepsilon }(\varepsilon +2)\text{csch}\left(\frac{\pi }{\varepsilon }\right)}{(\varepsilon +1)^2} \left(\sqrt{\varepsilon } \left(\sinh \left(\frac{\pi -t}{\sqrt{\varepsilon }}\right)-\sinh \left(\frac{\pi }{\sqrt{\varepsilon }}\right) \cos (t)\right)+\cosh \left(\frac{\pi }{\sqrt{\varepsilon }}\right) \sin (t)\right),
\end{multline}
while the optimal objective $J_*$ is found to be
\begin{equation}
J_*=\sqrt{\varepsilon } \tanh \left(\frac{\pi }{2 \sqrt{\varepsilon }}\right) \left(\frac{\coth ^2\left(\frac{\pi }{2 \sqrt{\varepsilon }}\right)}{(\varepsilon +1)^2}+1\right)-\frac{\pi  (2 \varepsilon +3)}{4 (\varepsilon +1)}.
\end{equation}

With these calculations in hand and further considering any  time $0<t\leq T$, we may take these optimal functions' limits as $\varepsilon$ approaches zero from above. We find that
\begin{subequations}\label{eq:limit}
\begin{align}
    \lim_{\varepsilon\to0^+}u_*&:=u_{\rm limit}= -1-\cos(t)\\
    \lim_{\varepsilon\to0^+}x_*&:=x_{\rm limit}= -\frac{1}{2} t \sin (t)+\cos (t)-1\,\\
    \lim_{\varepsilon\to0^+}J_*&:=J_{\rm limit}=-\frac{3\pi}{4}.
\end{align}
\end{subequations}
Although the unregularized problem exhibits no features of an ill-posed problem, note that the limiting optimal control fails to remain in the admissible class $\mathcal{U}$ since $u_{\rm limit}$ does not vanish at $t=0$. However, from the explicit calculation, we have that for each $\varepsilon>0,$ there exists a $u\in\mathcal{U}$ which solves the optimal control problem~\eqref{eq:testcontprob}. Thus, taking small values of the regularization parameter $\varepsilon$ manifests in a type of boundary layer in the vicinity of $t=0$ for optimal control $u_*$. That is, the optimal control develops a steep gradient at $t=0$ as $\varepsilon\to0^+,$ and, consequently, $\lim_{\varepsilon\to0^+}\norm*{u_*-u_{\rm limit}}_{H^1([0,T])}\neq0.$ The effect this has on numerical computations will be discussed shortly.

We numerically solve optimal control problem~\eqref{eq:testcontprob} in two ways. The first, performed pointwise in time, is to use MATLAB's \texttt{fmincon} which implements an interior-point method~\cite{BoydV}. The second is to use a method due to Calarco, et al.~\cite{Caneva,Doria}, called the $\mathbf{C}$hopped $\mathbf{Ra}$ndom $\mathbf{B}$asis (CRAB) method. It relies on choosing controls from the span of an appropriately chosen finite set of basis functions so that the optimization is performed over a set of unknown coefficients. The choice of basis is problem-specific but is always chosen so that controls remain in the appropriate admissible space. The choice of representation we make here is the following sine series:
\begin{equation}\label{eq:randclass}
    u_r(t)=\sum_{j=1}^{N}\frac{r_j}{j^2}\sin\left(\frac{2j\pi}{T}\right),\quad t\in[0,T],\\
\end{equation}
where the coefficients $r_j\in[-1,1]$ are parameters to be optimized over. Note that we choose the amplitudes in Equation~\eqref{eq:randclass} to decay quadratically because the Fourier series of absolutely continuous functions exhibits this type of decay~\cite{Trefethen}. In this way, the search space of amplitudes  is not severely restricted, yet early iterations of the controls generated by the 
CRAB method are not highly oscillatory.

\begin{figure}[htbp]
\begin{centering}
\subfigure{\includegraphics[width=.47\textwidth]{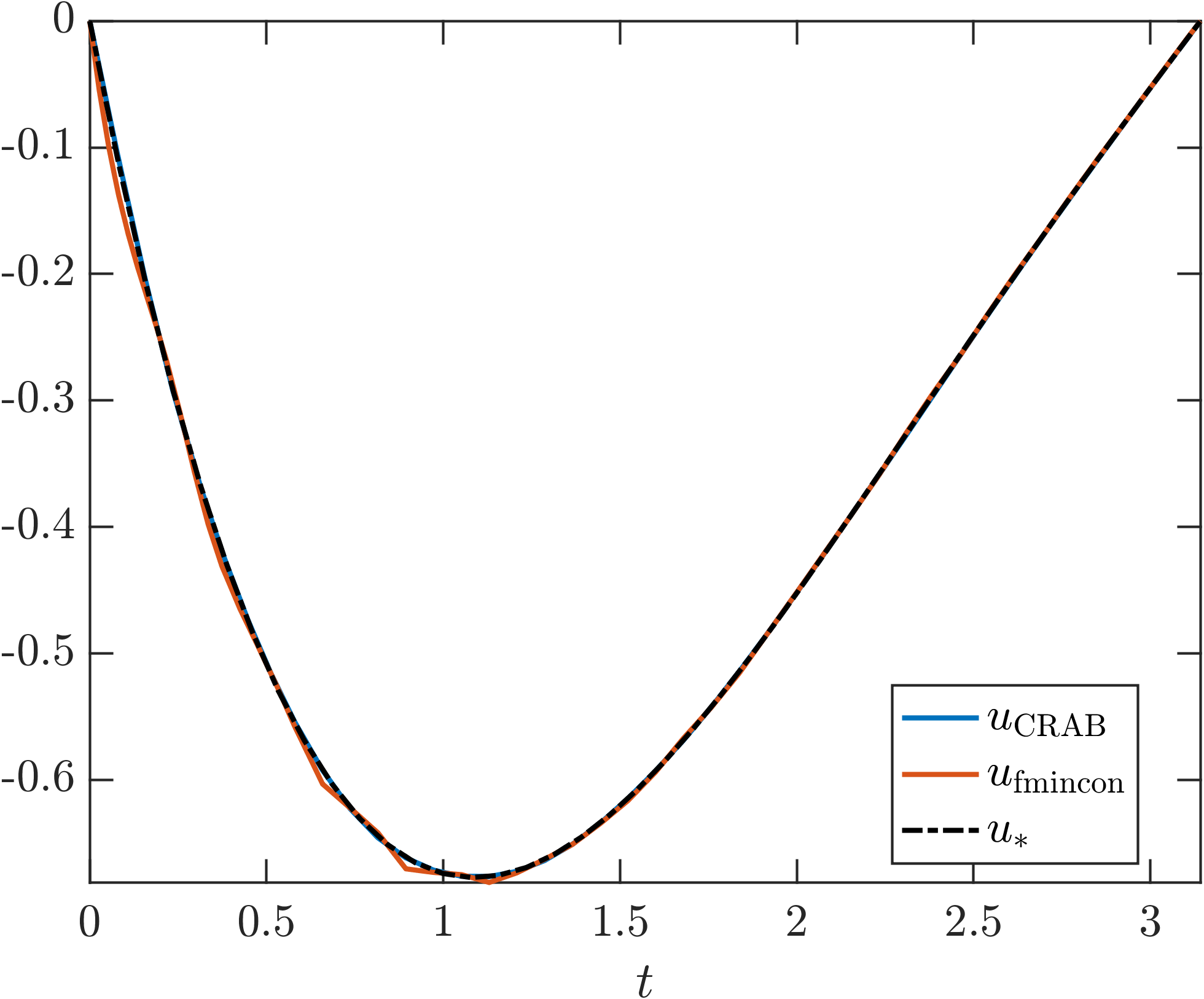}}
\subfigure{\includegraphics[width=.47\textwidth]{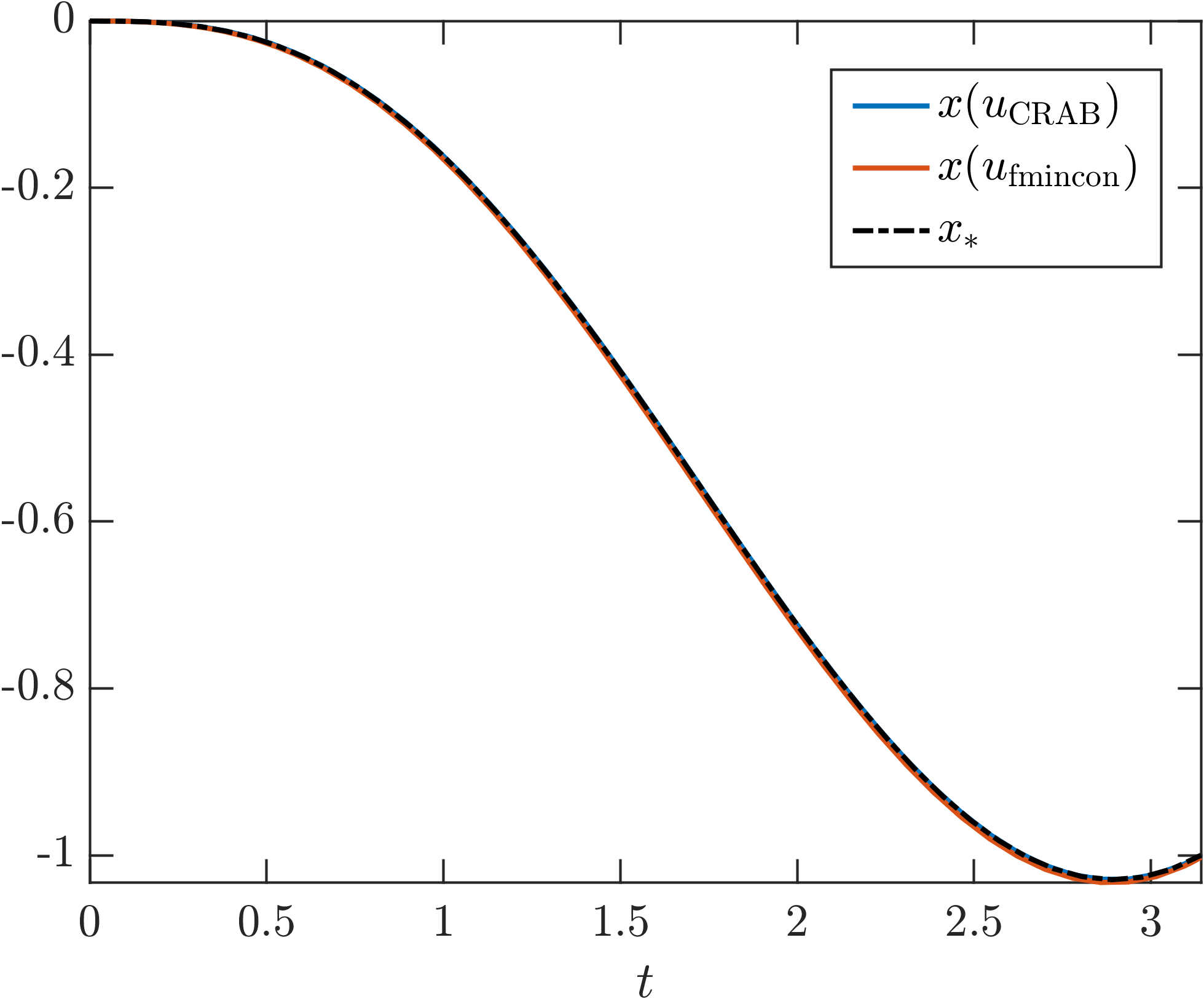}}
\caption{A numerical solution of optimal control problem~\eqref{eq:testcontprob}, with $\varepsilon$=1, via the CRAB method and MATLAB's \texttt{fmincon}. Left Panel: The controls found via both numerical methods and the exact optimal control given in closed-form by Equation~\eqref{eq:closedcontrol}. Right Panel: The states found via both numerical methods and the exact state given in closed-form by Equation~\eqref{eq:closedstate}.}
\label{fig:test1}
\end{centering}
\end{figure}

The CRAB method is of Galerkin-type. Recall that the intent behind any Galerkin method is to choose the number of basis functions $N$ simultaneously large enough to define an accurate approximation, yet small enough so that the overall procedure remains computationally inexpensive. For this problem, we find that just 12 modes are sufficient. The optimization problem is now a small-scale nonlinear programming (NLP) problem and can be solved using any number of industry-standard techniques. The technique we use to solve the resulting NLP problem is Differential Evolution (DE)~\cite{Storn}.

In our numerical study, we focus on investigating the role of the Tikhonov parameter $\varepsilon.$ To this end, we set $\varepsilon=1$ and show the numerical solution of Problem~\eqref{eq:testcontprob} in Figure~\ref{fig:test1}. We solve the state equations~\eqref{eq:states} with MATLAB's $\texttt{ODE45}$ solver. The numerics shown there demonstrate the adequacy of both methods used to solve Problem~\eqref{eq:testcontprob}. The numerical implementation of these methods is provided in the form of MATLAB code in Appendix~\ref{section:appendix}.

\begin{figure}[!ht]
\begin{centering}
\subfigure{\includegraphics[width=.47\textwidth]{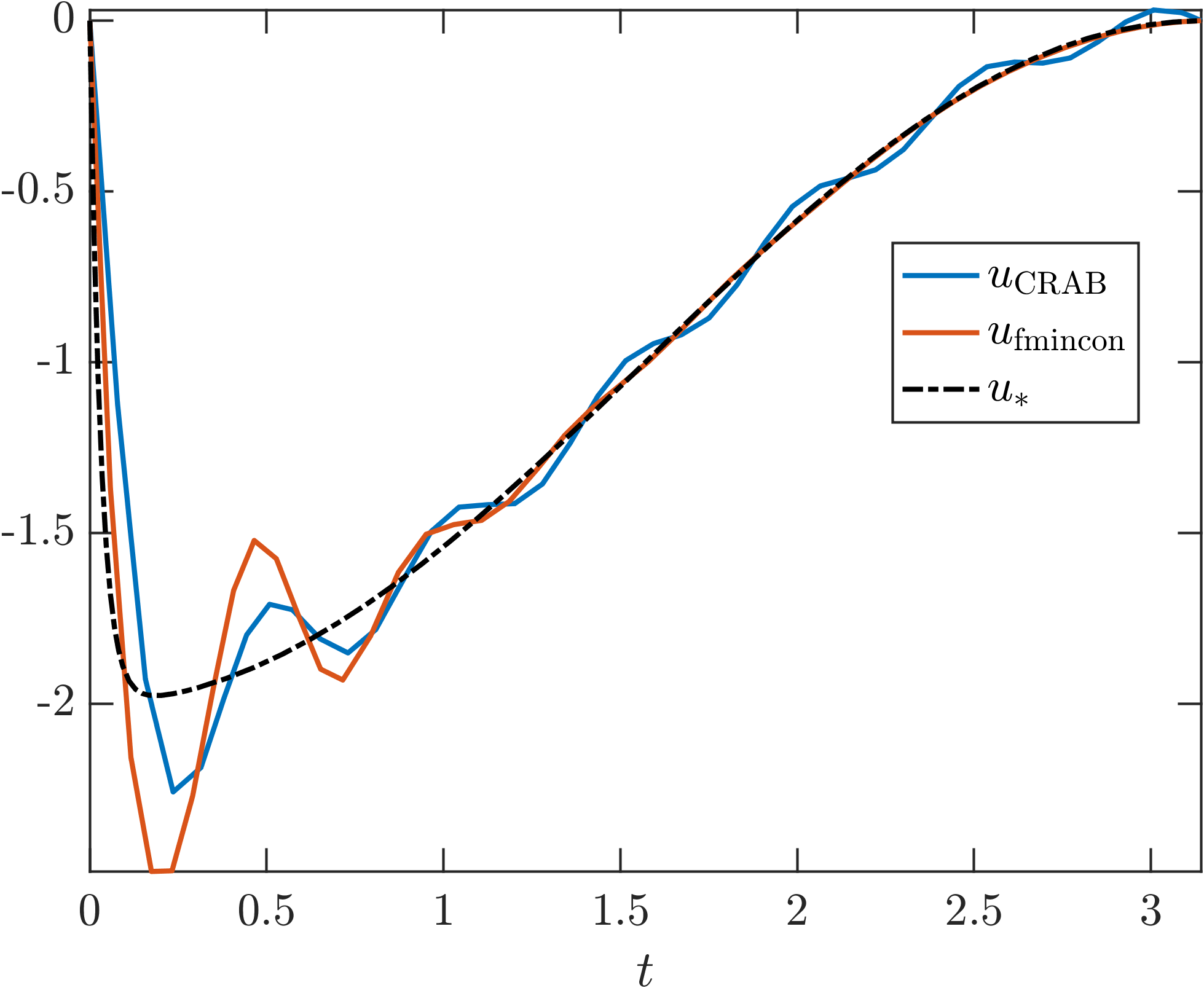}}
\subfigure{\includegraphics[width=.47\textwidth]{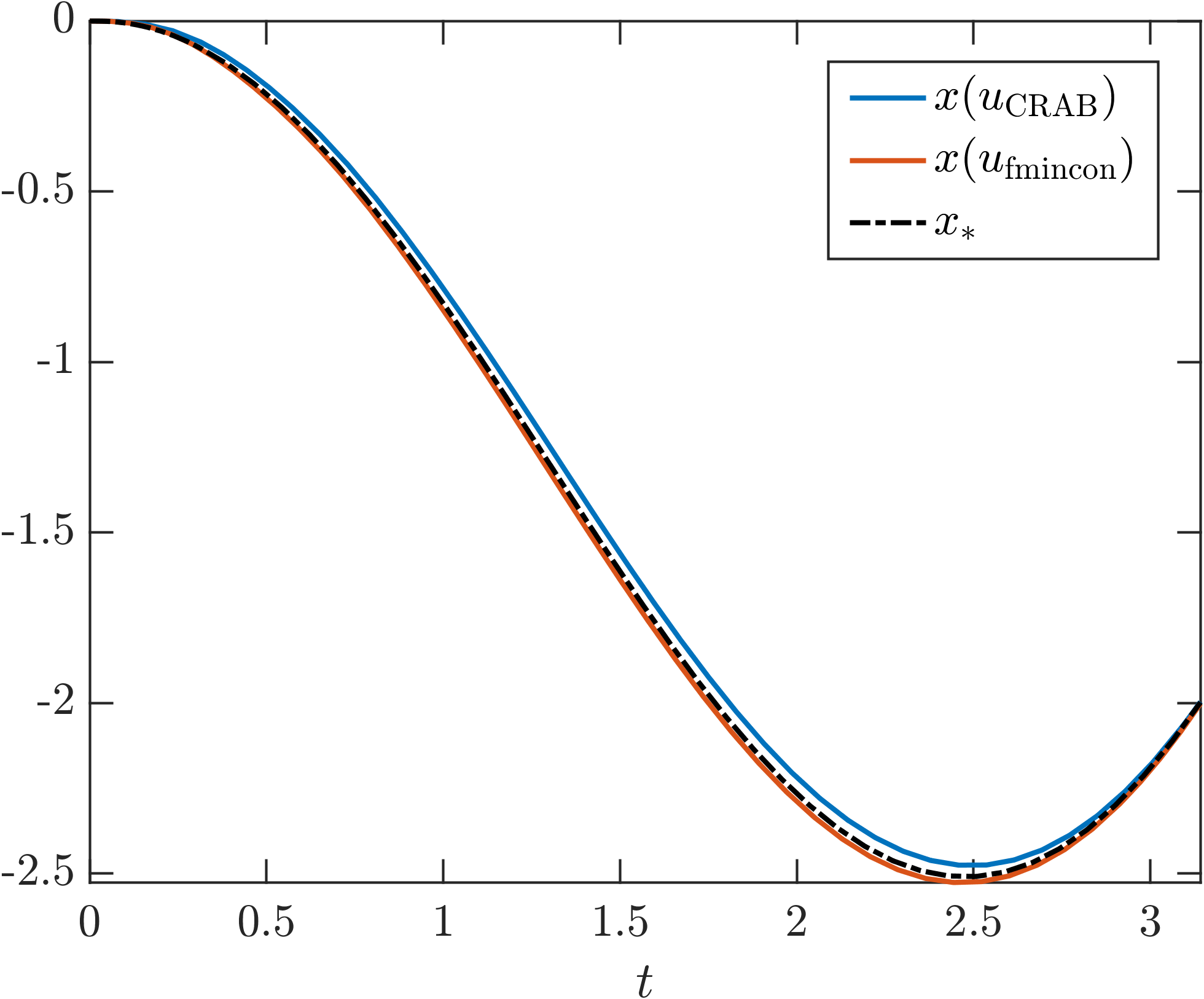}}
\subfigure{\includegraphics[width=.47\textwidth]{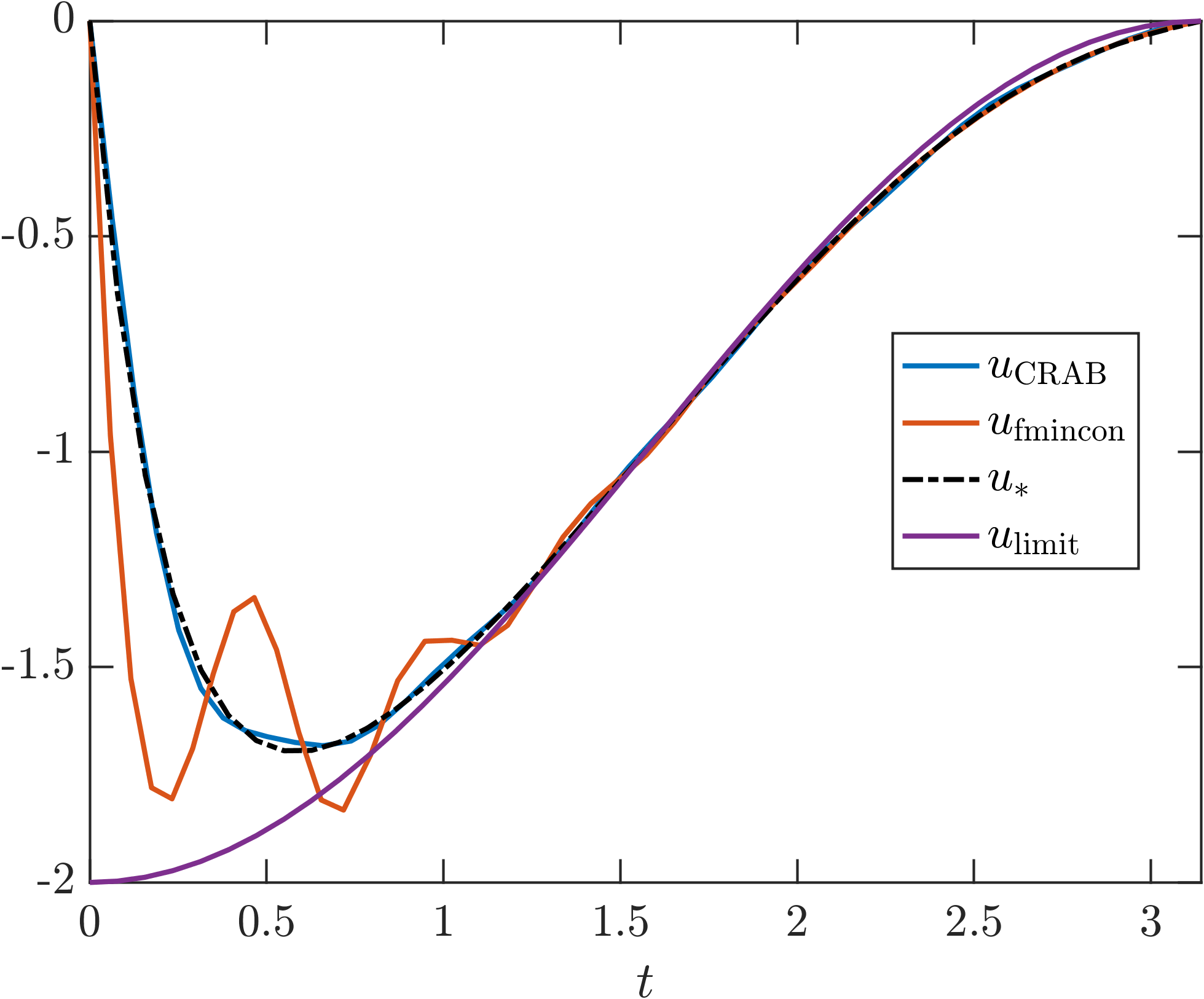}}
\subfigure{\includegraphics[width=.47\textwidth]{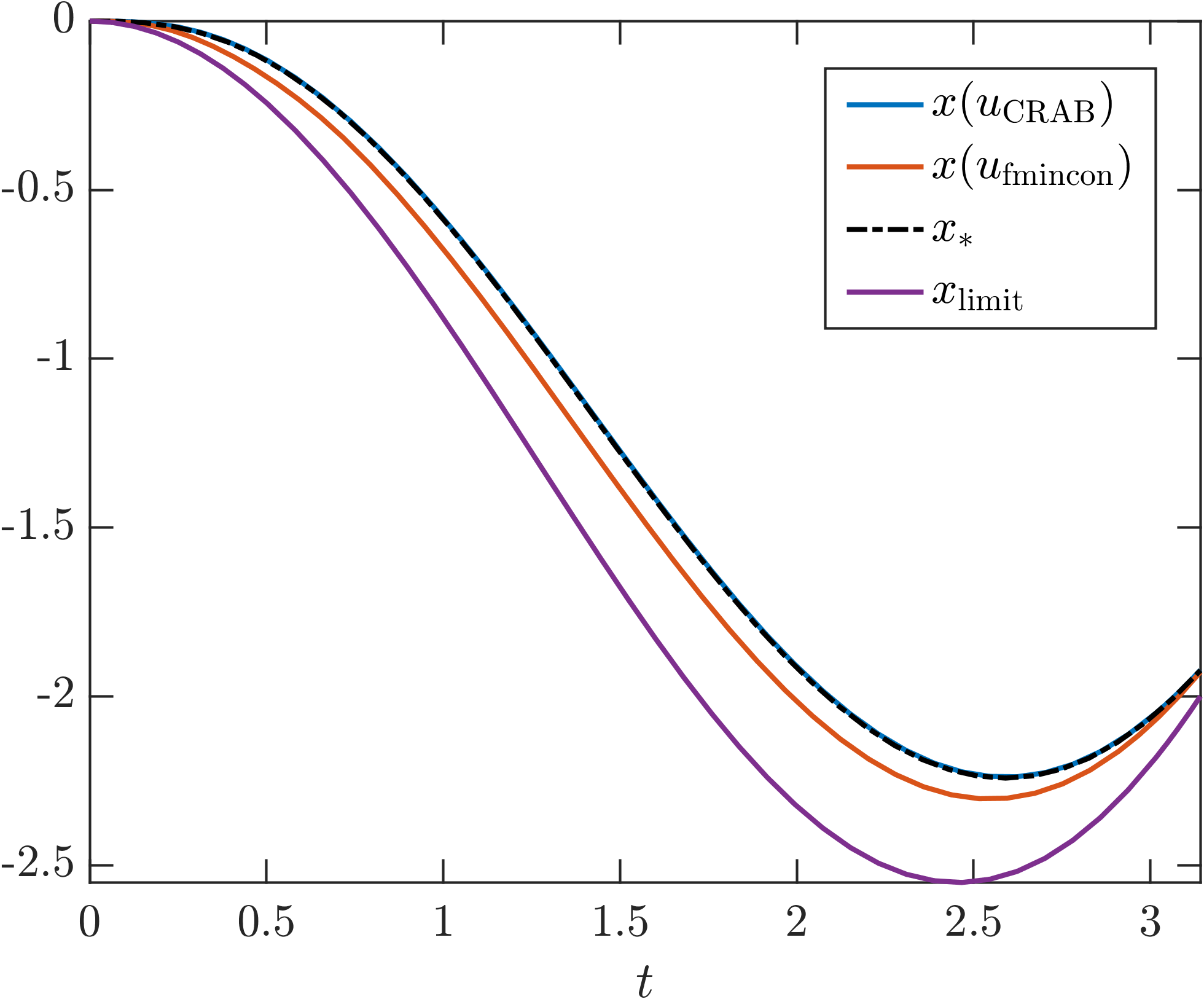}}
\caption{Examples of when the Tikhonov parameter $\varepsilon$ is small. The conventions here are consistent with Figure~\ref{fig:test1}. The top panels correspond to $\varepsilon=0.001$ while the bottom panels have $\varepsilon=0.04$. The bottom panels also include the limiting control and state functions given by Equation~\eqref{eq:limit}. }
\label{fig:test2}
\end{centering}
\end{figure}

By setting $\varepsilon=1,$ we find that we have strongly regularized the problem. This is seen in a pointwise comparison between $x_*,$ shown in Figure~\ref{fig:test1}, and $x_{\rm limit}$, shown in Figure~\ref{fig:test2}, especially for later times $t$. Yet, by taking $\varepsilon=0.001,$ we find that both optimization methods struggle to resolve the sharp transition of the control in the boundary layer near $t=0.$ We find $\varepsilon=0.04$ to be a good intermediate value under the tension of under- and over-regularizing the control problem (the bias-variance dilemma). In this case, we see that the CRAB method successfully finds the optimal control $u_*,$ while \texttt{fmincon} develops undesirable oscillations in the control near the boundary layer.

\section{Regularizing an Optical Beam Reshaping Problem}\label{section:Beam}
\subsection{Mathematical Formulation}
In this section, we solve an optimal control problem with an objective functional, first used in the context of high-fidelity quantum fluid manipulations by Hohenester, et al.~\cite{Hohenester}, here constrained by Schrödinger's equation in dimensionless form
\begin{equation}\label{eq:Intro2}
    i\psi_z=-\frac{1}{2}\triangle\psi+V(x,z)\psi.
\end{equation}
In applications outside of quantum mechanics, Schrödinger's equation is used as a standard model for paraxial light beam propagation as it arises from a slowly varying amplitude approximation to the variable-coefficient Helmholtz equation~\cite{JWGood}. Here, $z$ is the axis of propagation, $x$ is the transverse direction, $\triangle$ is the Laplacian in $x$, $V(x,z)$ is proportional to the square of a spatially varying refractive index, and the wavefunction $\psi(x,z)$ is interpreted as a spatially varying complex electric amplitude.  We assume the propagation media is lossless, hence the potential $V$ is a real function of the waveguide coordinates. 

We introduce a parametrization $w(z)$ for the potential $V$ in order to reduce the dimensionality of the numerical control problem for ease of computation. Thus, the light reshaping problem we consider is as follows: Find the optimal control $w(z)$ that best transforms the intensity distribution of an initial Schrödinger state $\varphi_0(x)$
into the intensity distribution of the desired state $\varphi_d(x)$ satisfying
\begin{subequations}\label{eq:initdes}
\begin{align}
-\frac{1}{2}\triangle\varphi_0(x)+V(x,w(0))\varphi_0(x)&=\lambda_0\varphi_0(x),\\
-\frac{1}{2}\triangle\varphi_d(x)+V(x,w(l))\varphi_d(x)&=\lambda_l\varphi_d(x),
\end{align}
\end{subequations}
i.e., the initial and desired states are eigenfunctions, of the time-independent Schrödinger operator $P=-\frac{1}{2}\triangle+V(x,w(z))$, at $z=0$ and at the end of a specified propagation length $l$, respectively. Thus, we formulate the problem of designing an optimal coupler between two waveguides with different transverse profiles and their eigenpairs $\left(\varphi_0,\lambda_0\right)$ and $\left(\varphi_d,\lambda_l\right)$. 

The following problem structure is due to Hohenester, et al.~\cite{Hohenester}, and uses the following objective functional
\begin{equation}\label{eq:BeamObj}
J=\frac{1}{2}\left(\|\varphi_d(\cdot)\|_{L^2(\mathbb{C}^n)}^4-\left|\left\langle \varphi_d(\cdot),\psi(\cdot,l)\right\rangle_{L^2(\mathbb{C}^n)}\right|^2\right)+\frac{\varepsilon}{2}\int_0^l\left|\partial_z w\right|^2dz,
\end{equation}
where $\varepsilon>0$ and $z\in(0,l)$ with $l>0$. The objective functional $J$ involves the \emph{infidelity} 
\begin{equation}\label{eq:Jinfed}
    J_{\rm infidelity}= \frac{1}{2}\left(\|\varphi_d(\cdot)\|_{L^2(\mathbb{C}^n)}^4-\left|\left\langle \varphi_d(\cdot),\psi(\cdot,l)\right\rangle_{L^2(\mathbb{C}^n)}
    \right|^2\right)
\end{equation} 
which penalizes misalignments of the computed function $\psi(x,l)$ with respect to the desired state $\varphi_d(x)$. In the language of optimal control theory~\cite{Bryson,Calculus1989}, the infidelity is called a terminal cost. This objective functional disregards the physically unimportant global phase difference between the desired and computed states,  a significant advantage over a typical least-squares approach. The second contribution to the objective, the running cost over $[0,l]$, is the Tikhonov regularization, analogous to the regularization in Objective~\eqref{eq:testcontprob}. Hintermuller, et al., prove the control framework of Hohenester, et al., is well-posed with the introduction of this Tikhonov regularization, i.e., there exists a control $w\in H^1([0,l])$ that minimizes the objective $J$~\cite{Hintermuller}.

The search for optimal controls is performed over the admissible class 
\begin{equation}
    \mathcal{W}=\left\{w\in H^1\left([0,l]\right):w(0)=w_0,w(l)=w_l\right\}.
\end{equation}
We assume the eigenfunctions $\varphi_d$ and $\varphi_0$ are both in the space $H^{1}(\mathbb{R}^n)$. We also assume that the eigenfunctions $\varphi_0$ and $\varphi_d$ have unit intensity, i.e., $\|\varphi_0\|_{L^2({\mathbb{R}^n})}=\|\varphi_d\|_{L^2({\mathbb{R}^n})}=1$, so that the infimum of the infidelity~\eqref{eq:Jinfed} is 0. Lastly, we assume the potential $V(x,u(z))$ is in the space of bounded, continuous functions, i.e., $C_b^{0}([0,l];H^{1}\left(\mathbb{C}^n\right)),$ for every $w\in\mathcal{W}$. With the above assumptions in place, the regularity of the wavefunction $\psi$ solving Equation~\eqref{eq:Intro2} is known~\cite{MASPER}; $\psi\in C^1([0,l];H^1(\mathbb{C}^n))$. Moreover, the control problem with objective functional~\eqref{eq:BeamObj} is well-posed for sufficiently large $\varepsilon>0$~\cite{Hintermuller}. 

By using the method of Lagrange multipliers and the fundamental theorem of calculus, the optimal control problem can be written as
\begin{equation}\label{eq:Lastbeamprob}
\min_{w\in \mathcal{W}}J=
\min_{w\in \mathcal{W}}\int_0^l\mathcal{L}(\psi,\partial_z\psi,\triangle{\psi},\psi^{\dag},p^{\dag},w,\partial_zw) dz,
\end{equation}
where the Lagrange density is given by
\begin{equation}
\begin{split}
    \mathcal{L}&=\Re\left\{\left<p,i\partial_z\psi+\frac{1}{2}\triangle\psi-V(x,w(z))\psi\right>-\left<\varphi_d,\psi\right>\left<\partial_z\psi,\varphi_d\right>\right\}+\frac{\varepsilon}{2}|\partial_zw|^2,\\
\end{split}
\end{equation}
where $p$ is a Lagrange multiplier with inner products understood in the sense of $L^2(\mathbb{C}^n)$ and daggers denoting the Hermite conjugate. Again, it is straightforward to show, using standard arguments from the calculus of variations, that the optimality conditions of Problem~\eqref{eq:Lastbeamprob} are given by
\begin{subequations}\label{eq:optcond}
\begin{alignat}{2}
\label{eq:stateeq}
i\partial_z\psi&=-\frac{1}{2}\triangle\psi+V(x,w(z))\psi,\qquad &\psi(x,0)=\varphi_0(x), \\
\label{eq:costateeq}
i\partial_zp&=-\frac{1}{2}\triangle p+V(x,w(z)) p,\qquad &ip(x,l)=\left\langle \varphi_d,\psi(x,l)\right\rangle_{L^2(\mathbb{C}^n)}\varphi_d, \\
\label{eq:controleq}
\varepsilon\partial_z^2{w}&=-\Re\left\langle p,\partial_wV\psi\right\rangle_{L^2(\mathbb{C}^n)},\qquad &w(0)=w_0,\ w(l)=w_l.
\end{alignat}
\end{subequations}
Equation~\eqref{eq:costateeq} is the adjoint equation of Equation~\eqref{eq:stateeq} and governs the axial evolution of the costate $p(x,t)$ backward from its terminal condition at $z=l$. 

To state the well-posedness of the control problem more precisely, consider the so-called reduced objective functional
\begin{equation}\label{eq:reducedcost}
\mathcal{J}:\mathcal{W}\to\mathbb{R},\quad w\mapsto\mathcal{J}[w]:=J\left[\psi(w),w\right].
\end{equation}
Let $w_*$ denote the optimal control, and define $\psi_*:=\psi(w_*)$, $p_*:=p(w_*)$. Since the optimal control problem is well-posed, then for every $w\in\mathcal{W}$,
\begin{equation}\label{eq:optineq}
    \mathcal{J}[w]\geq \mathcal{J}[w_*]=\min_{w\in\mathcal{W}}\mathcal{J},
\end{equation}
i.e., the minimum is attained by the optimal control $w_*$. In addition, at the minimum, the optimal triple $(\psi_*,p_*,w_*)$ satisfies
Equations~\eqref{eq:optcond}. For this reason, pursuing numerical approximations of Equations~\eqref{eq:optcond} and the optimality condition~\eqref{eq:optineq} when searching for the optimal control $w_*$ is meaningful.

The similarity of Equation~\eqref{eq:stateeq} and Equation~\eqref{eq:costateeq} is due to the self-adjoint nature of the Schrödinger operator $P= -\frac{1}{2}\triangle +V(x,w(z))$.  Equation~\eqref{eq:controleq} governs the optimal control $w$, and together with the boundary conditions defined through the admissible class $\mathcal{W}$, is a boundary value problem on $[0,l]$. We solve Equations~\eqref{eq:stateeq} and~\eqref{eq:costateeq} via a second-order Fourier split-step method~\cite{Strang}, where the $z-$dependence of the potential is handled by the midpoint method. Equation~\eqref{eq:controleq} is not solved numerically but is instead reinterpreted in the context of a projected gradient descent method that ensures controls remain admissible~\cite{me,me2,Borzi}.

\subsection{Numerical Example}
In order to demonstrate the role Tikhonov regularization has on light reshaping in a simple setting, we set $x\in\mathbb{R},$ i.e. $n=1,$ and consider initial and terminal eigenfunctions for which $V(x,0)$ and $V(x,l)$ can be computed in closed form. It is well-known that the so-called Pöschl-Teller potential, 
\begin{equation*}
    V(x)=-\frac{s(s+1)}{2}\sech^2(x),
\end{equation*}
$s\in\mathbb{N},$ gives Legendre functions as eigenfunctions for the time-independent Schrödinger equation\cite{Poschl}. As a test, we consider the problem of reshaping the ground state eigenfunction for $s=1$ to the ground state corresponding to $s=4.$ We introduce two parametrizations for the reshaping potential; a "depth" control $u(z)$ and a "width" control $v(z)$. More precisely, we assume the following form of the potential 
\begin{equation}\label{eq:PosclPot}
V(x,u(z),v(z))=-\frac{u(z)}{2}\sech^2(v(z)x),
\end{equation}
where the initial and terminal eigenfunctions are given by
\begin{equation}\label{eq:poscleig}
\varphi_0(x)=-\frac{1}{\sqrt{2}}\sech(x),\quad \varphi_d(x)=-\frac{3}{2\sqrt{3}}\sech^2(x),
\end{equation} 
and the appropriate control boundary conditions are
\begin{subequations}
\begin{align}
\label{eq:uBC}
&u(0)=2,\quad u(l)=20,\\
\label{eq:vBC}
&v(0)=1,\quad v(l)=1.
\end{align}
\end{subequations} 
This assumption on $V(x,z)$ slightly changes the optimality condition~\eqref{eq:controleq} such that the following equations
\begin{subequations}
\begin{alignat}{2}\label{eq:controlueq}
\varepsilon\partial_z^2{u}&=-\Re\left\langle p,\partial_uV\psi\right\rangle_{L^2(\mathbb{C})},\qquad &u(0)=u_0,\ u(l)=u_l,\\
\varepsilon\partial_z^2{v}&=-\Re\left\langle p,\partial_vV\psi\right\rangle_{L^2(\mathbb{C})},\qquad &v(0)=v_0,\ v(l)=v_l,
\end{alignat}
\end{subequations}
are now the appropriate Euler-Lagrange equations for the controls $u$ and $v$, while the state and costate equations~\eqref{eq:stateeq},~\eqref{eq:costateeq} remain unchanged.
\begin{figure}[!ht]
\begin{centering}
\subfigure{\includegraphics[width=.47\textwidth]{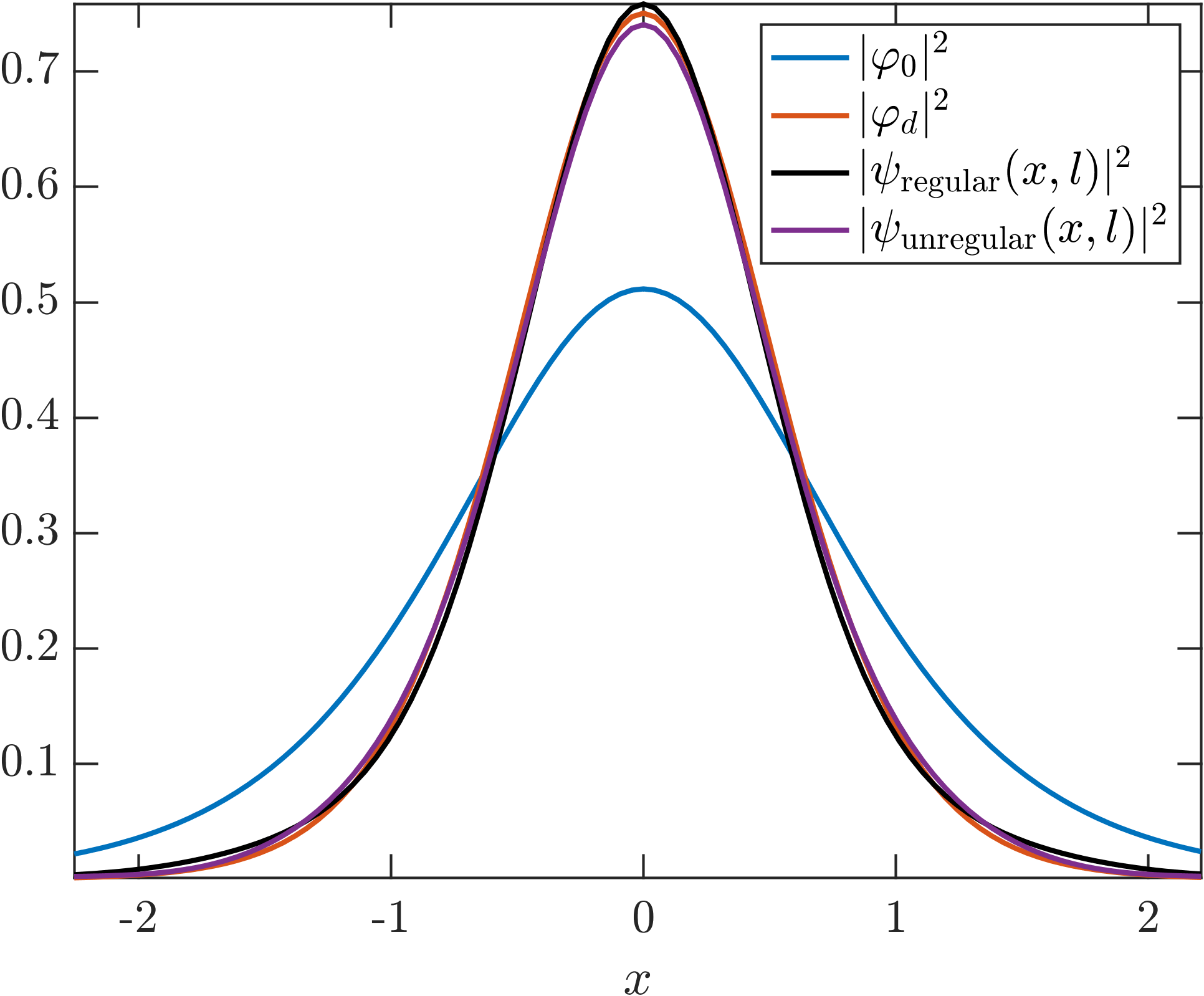}}
\subfigure{\includegraphics[width=.47\textwidth]{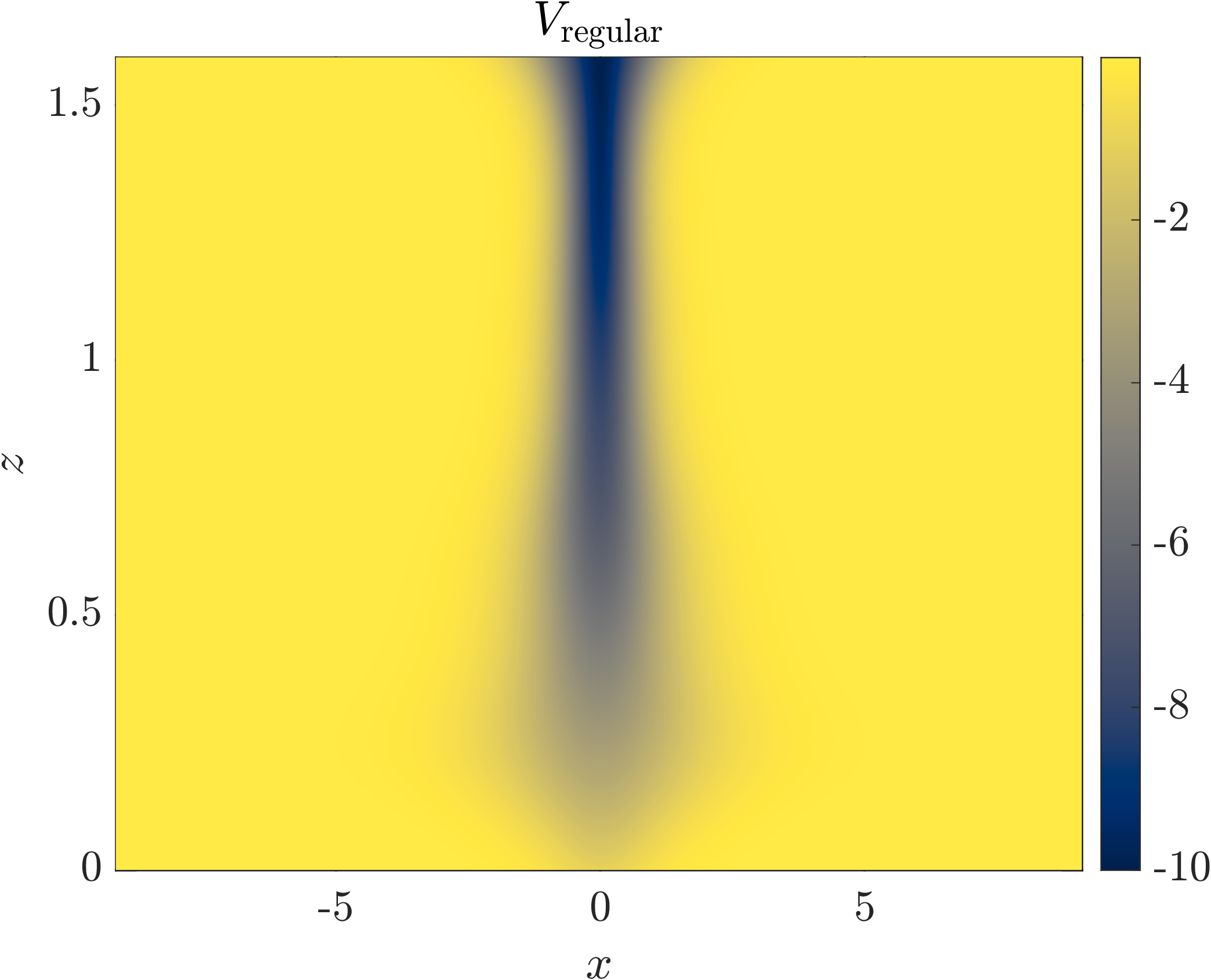}}
\subfigure{\includegraphics[width=.47\textwidth]{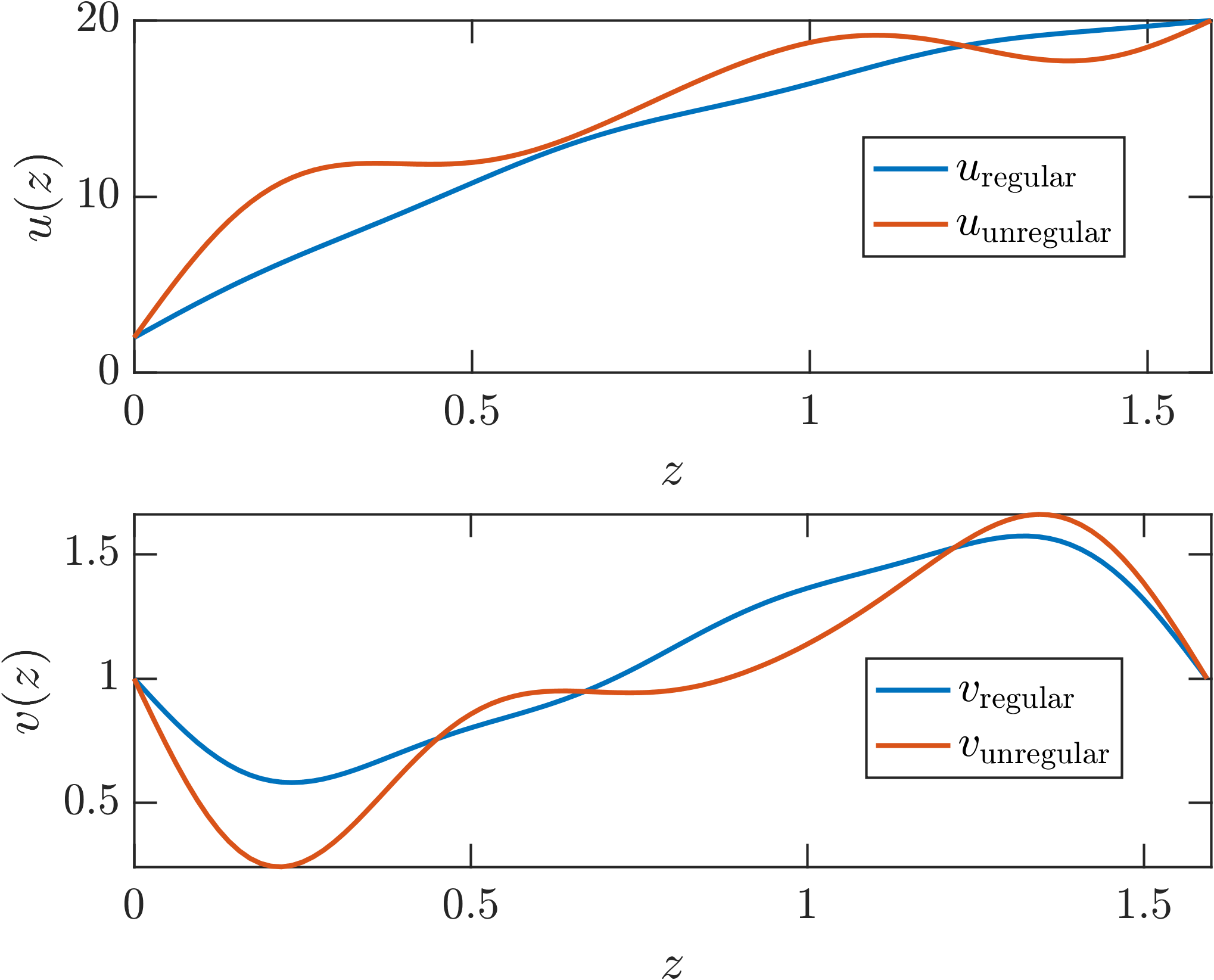}}
\subfigure{\includegraphics[width=.47\textwidth]{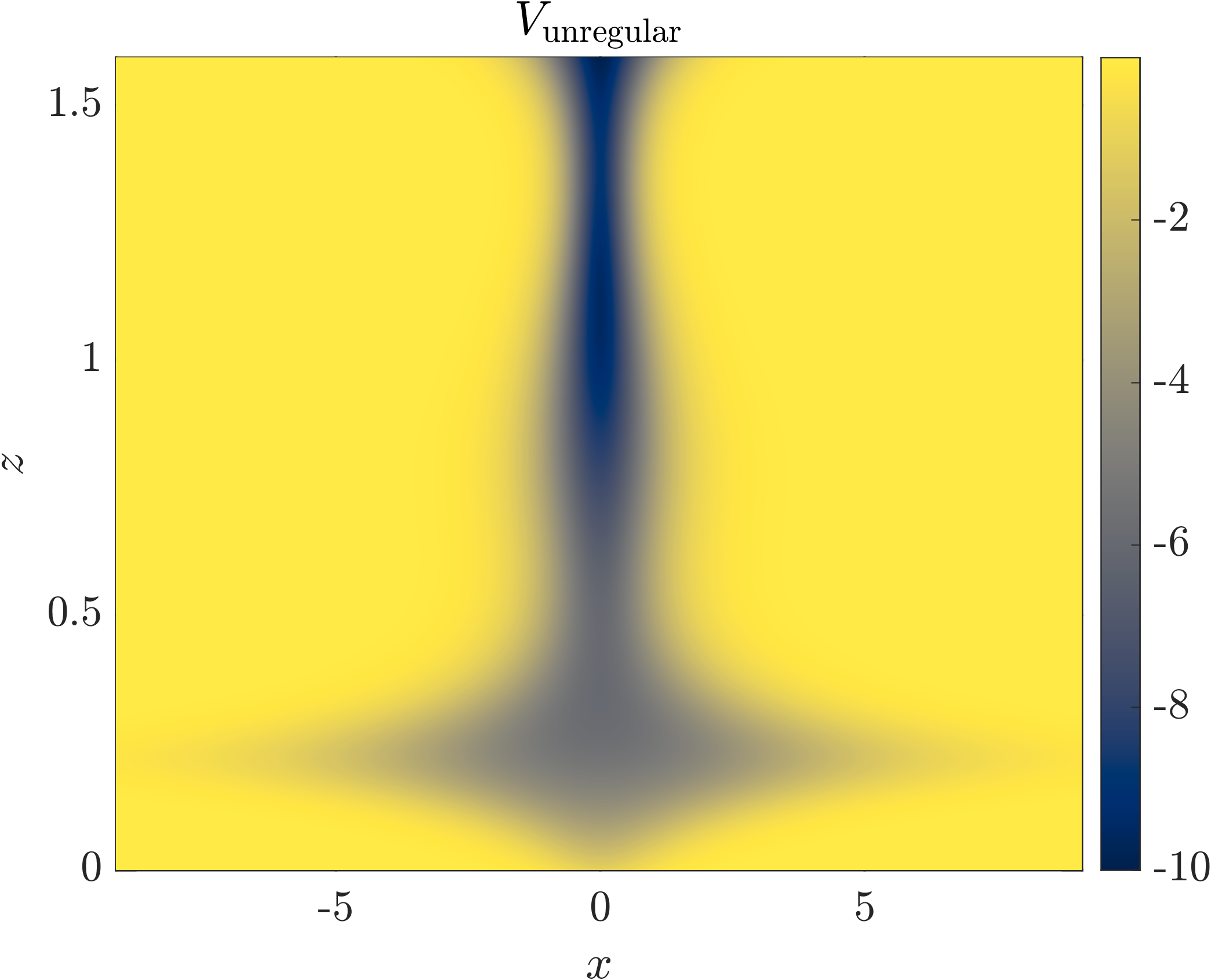}}
\caption{Tikhonov regularization applied to light reshaping with Objective~\eqref{eq:BeamObj} and a regularization parameter of  $\varepsilon=10^{-4}.$ Top left: intensity distributions of the initial, desired, and computed wavefunctions. Bottom left: numerically computed regularized and unregularized optimal controls. Right column: the resulting regularized and unregularized potentials $V.$ We see that the regularized potential yields a satisfactory result in the reshaping of the intensity distribution of light, yet is more feasible to fabricate than the unregularized potential.}
\label{fig:Schrod}
\end{centering}
\end{figure}

For numerical implementation, we choose the following domain and discretization parameters: $l=1.6,$ $l\times2^7$ discretization points in $z,$ and $x\in[-4\pi,4\pi]$ with $2^{11}$ discretization points.  We choose $l=1.6$ because this is the smallest length found by experimentation such that the infidelity is less than $10^{-3}$ for the unregularized case of when $\varepsilon=0$. Having a smaller length $l$ in the design is, of course, desirable since such designs are more cost-effective in manufacturing.

In Figure~\ref{fig:Schrod}, we show the result of using the CRAB method, with 10 basis functions for each control $u$ and $v$ about binomials satisfying the appropriate boundary conditions on $u$ and $v$, coupled with a projected gradient descent method~\cite{me2} in order to solve the optimal control problem. We find that introducing a small Tikhonov regularization weight of $\varepsilon=10^{-4}$ greatly impacts the control problem. Despite the infidelities of the regularized and unregularized control problems remaining well within the same order of magnitude, the regularized potential $V(x,u(z),v(z))$ has tighter support about the origin throughout the length of the optical device. We see in this case that Tikhonov regularizations aid in discovering designs that are more cost-effective in fabrication.

\section{Concluding Remarks}
In this work, we provide standard optimization problems where applying the Tikhonov regularization technique is crucial. Moreover, the first couple is instructive and can be used in introductory courses on numerical linear algebra and optimization. We also demonstrate the utility of the CRAB method, over a more standard numerical optimization method, when trying to solve control problems that exhibit sharp transition regions in the optimal control. Finally, we demonstrate the utility of Tikhonov regularization for a more realistic design problem. We see that regularization does well to smooth the fast transitions exhibited by the controls near their boundary. This boundary layer-like phenomenon may be difficult to understand in this more complex setting, yet is nicely illustrated in the second optimization problem given in this paper.

\section*{Acknowledgments}
The author gratefully acknowledges helpful discussions with M. Siegel, M. Porter, and  T. Witelski for their review and feedback of the manuscript. The author especially acknowledges feedback and careful review of the work which resulted in this manuscript from R.H. Goodman while the author was still a student at the New Jersey Institute of Technology.

\section*{Disclosures} The authors declare no conflicts of interest.

\section*{Data Availability} Data underlying the results presented in this paper are not publicly available at this time but may be obtained from the authors upon reasonable request.

%%%%%%%%%%%%%%%%
% Bibliography %
%%%%%%%%%%%%%%%%
\bibliographystyle{abbrv} % We choose the "plain" reference style
\bibliography{bib.bib} % Entries are in the "bibliography.bib" file
\newpage
\appendix
%\definecolor{bg}{rgb}{0.95,0.95,0.95}
%[linenos=true,bgcolor=bg]{matlab}
\section{Main Script Accompanying Section~\ref{section:OC}}\label{section:appendix}
\begin{lstlisting}[language=Matlab]
%% Solves an Optimal Control Problem via CRAB and fmincon
clear, close all

% Tikhonov weight g and Optimal control in closed-form
g=.04;
uExact=@(x)(-1).*(1+g).^(-1).*(1+g+cos(x)+(-1).*csch(g.^(-1/2).*pi).*((2+g).* ...
  sinh(g.^(-1/2).*(pi+(-1).*x))+g.*sinh(g.^(-1/2).*x)));

% Initial conditions for the state and control boundary conditions
x0=[0 0];
u0=0;uT=0;t0=0;T=pi;

%Define the objectives for each case and number of optimization points
JCRAB=@(u)getObjective(u,g,false,x0,t0,T);
Jfmin=@(u)getObjective(u,g,true,x0,t0,T);

% Differential Evolution Parameters
DEParams.F=0.8;DEParams.CR=0.7;
DEParams.NP=120;DEParams.Nmax=300;
DEParams.ND=12;

% Amplitude hyperrectangle for the CRAB method
l=zeros(1,DEParams.ND);
for j=1:DEParams.ND
    l(j)=(2*rand-1)/j^2;
end

% Execute differential evolution
optimalamps=diffevoND(JCRAB,l,DEParams);

% Reconstruct the control from the CRAB amplitudes and compute resulting state
uCRAB=recon(linspace(t0,T,1000),optimalamps(:,end));

% Prepare for a blackbox optimization
options = optimoptions(@fmincon,'Display','iter','Algorithm','interior-point');

% Keep number of function evaluations consistent with CRAB
options.MaxFunctionEvaluations=DEParams.NP*DEParams.Nmax;
fminPts=15;
t=linspace(t0,T,fminPts);

% initial guess on the control
u00=t.*(t-T);

% set equality constraints on the control
Aeq = zeros(fminPts);
Aeq(1,1)=1;
Aeq(end,end)=1;
beq=zeros(fminPts,1);

% Execute an interior-point method using MATLAB's fmincon
ufmin=fmincon(Jfmin,u00,[],[],Aeq,beq,[],[],[],options);
ufmin=spline(t,ufmin);
ufmin=@(t)ppval(t,ufmin);

% Solve for the states
xCRAB =SolveStates(t0,T,uCRAB,x0);
xfmin =SolveStates(t0,T,ufmin,x0);
xExact=SolveStates(t0,T,uExact,x0);

\end{lstlisting}

\section{Functions Used}

%\definecolor{bg}{rgb}{0.95,0.95,0.95}
%[linenos=true,bgcolor=bg]{matlab}
\begin{lstlisting}[language=Matlab]
function x=SolveStates(t0,T,u,x0)
% Solve for states using ODE45
[~,states]=ode45(@(t,x)StateRHS(t,x,u),[t0 T],x0);
x=states(:,1);x=spline(t,x);x=@(t)ppval(t,x);
end
\end{lstlisting}
\begin{lstlisting}[language=Matlab]
function dxdt=StateRHS(t,x,u)
dxdt=[x(2)
    u(t)-x(1)];
end
\end{lstlisting}
\begin{lstlisting}[language=Matlab]
function J=getObjective(u,g,blackbox,x0,t0,T)
if blackbox
    % this is coming from fmincon
    t=linspace(t0,T,length(u));
    dt=t(2)-t(1);
    u=spline(t,u);
    u=@(t)ppval(t,u);
else
    % this is coming from CRAB
    t=linspace(t0,T,1000);
    dt=t(2)-t(1);
    u=recon(t,u);
end
% solve for the states
x=SolveStates(t0,T,u,x0);
% compute the derivative of the control
udot=gradient(u(t),dt);
% the objective functional
J=trapz(t,x(t)+.5*u(t).^2+.5*g*udot.^2);
end
\end{lstlisting}
\begin{lstlisting}[language=Matlab]
function u=recon(t,params)
% reconstruct the control from the CRAB amplitudes
t0=t(1);T=t(end)-t0;
sinepart=0;
for i=1:length(params)
    sinepart=sinepart+params(i)*sin(i*pi*(t-t0)'/T);
end
u=sinepart;
u=spline(t,u);
u=@(t)ppval(t,u);
end
\end{lstlisting}
\begin{lstlisting}[language=Matlab]
function optimalamps=diffevoND(f,l,DEparams)
%Minimizes f(x) using the differential evolution algorithm
%of Storn, Price
%Inputs:
%f is the function handle we want to minimize
%l is the length of the interval for the defined function
%DEparams.CR,DEparams.NP,DEparams.F are the diffevo
%parameters discussed in the literature
%Output: the possible global minimizer
%Initialize a population of candidate minimizers from the intervals [-l,l]
pop=zeros(DEparams.ND,DEparams.NP);
for j=1:DEparams.ND
    for i=1:DEparams.NP
        pop(j,i)=(l(j)*(2*rand-1));
    end
end
fpop=zeros(DEparams.NP,1);
z=zeros(DEparams.ND,1);
for i=1:DEparams.NP
    fpop(i)=f(pop(:,i));
end

%Set DEparams.tolerance and maximum number of iterations
counter=0;tic
while counter<DEparams.Nmax
    newpop=pop;
    newfpop=fpop;
    for i=1:DEparams.NP
        for j=1:DEparams.ND
            currentmember=pop(j,i);
            %Choose 3 distinct elements different from current member
            randind=randsample(DEparams.NP,3);
            while length(unique(randind))<3||~isempty(find(randind==i,1))
                randind=randsample(DEparams.NP,3);
            end
            a=pop(j,randind(1));
            b=pop(j,randind(2));
            c=pop(j,randind(3));
            
            %Compute currentmember's potentially new location y
            r=rand;
            %Potential Cross Over
            if r<DEparams.CR
                y=a+DEparams.F*(b-c);
            else
                y=currentmember;
            end
            z(j)=y;
        end
        fCandidate=f(z);
        if fCandidate<fpop(i)
            newpop(:,i)=z;
            newfpop(i)=fCandidate;
        end
    end
    %Update
    counter=counter+1;
    pop=newpop;
    fpop=newfpop;
    fprintf('%f cost after %i iterations of DE complete in %f seconds.\n',...
        min(fpop),counter,toc)
    [~,optindex]=min(fpop);
    optimalamps(:,counter)=pop(:,optindex);
end
end

\end{lstlisting}

\end{document}